# Addition Formulas of Leaf Functions According to Integral Root of Polynomial Based on Analogies of Inverse Trigonometric Functions and Inverse Lemniscate Functions

Kazunori Shinohara

Department of Mechanical System Engineering, Daido University
10-3 Takiharu-cho, Minami-ku, Nagoya 457-8530, Japan

**Abstract**

Inverse trigonometric functions *arcsin(r)* and *arccos(r)* are obtained by integrating the reciprocals of the root of the polynomial $1 - t^2$. Based on this analogy, the inverse lemniscate functions *arcsl(r)* and *arccl(r)* are created by integrating the reciprocal of the root of the polynomial $1 - t^4$. The lemniscate functions are then extended to the Jacobi elliptic functions, which in turn are improved using the theta functions. However, to the best of our knowledge, the integration of the reciprocals of the root of the polynomial $1 - t^6$ has never been published. The inverse functions based on $1 - t^6$ are defined as $arcsleaf_3(r)$ and $arccleaf_3(r)$, and they were determined to produce waves different from those of the trigonometric and lemniscate functions. This paper presents the addition formulas of the artificially created functions $sleaf_3(l)$ and $cleaf_3(l)$. These formulas were numerically verified through examples.

**Keywords**: Leaf functions, Trigonometric functions, Lemniscatic elliptic functions, Inverse functions, Addition formulas

## 1 Introduction

The second derivative of a function *r(l)* with respect to a variable *l* is equal to *-n*

times the function raised to the *2n-1* power of *r(l)* using this definition, an ordinary differential equation is constructed. Graphs with the horizontal axis as the variable *l* and the vertical axis as the variable *r(l)* are created by numerically solving the ordinary differential equation. These graphs show several regular waves with a specific periodicity and waveform depending on the natural number *n* [1] [2].

$$\frac{d^2 r(l)}{dl^2} = -n \cdot r(l)^{2n-1} \tag{1.1}$$

The definition of the arc leaf function $arcsleaf_n(r)$ is obtained by the above equation.

$$arcsleaf_n(r) = \int_0^r \frac{dt}{\sqrt{1-t^{2n}}} = l \tag{1.2}$$

In case of *n=3*, the above ordinary differential equations are as follows:

$$\frac{d^2 r(l)}{dl^2} = -3 \cdot r(l)^5 \tag{1.3}$$

These equations are related in the nonlinear equations (Quintic Duffing equation) [3]. To obtain the solution accuracy, or to grasp behavior of the solution, the addition formula often be demanded in the solution process [4] [5] [6]. Historically, the trigonomic function (in *n=1* of the Eq. (1.2)) and the lemniscatic elliptic function (in *n=2* of the Eq. (1.2)) have been discussed [7] [8] [9] [10] [11]. However, in case of *n=3*, the integration of the reciprocals of the root of the polynomial $1-t^6$ has never been published.

The variable *n* represents a natural number. For *n = 1*, the above equation is written as follows:

$$arcsleaf_1(r) = \int_0^r \frac{dt}{\sqrt{1-t^2}} = \arcsin(r) = l \tag{1.4}$$

The function *arcsin(r)* represents inverse trigonometric functions. The double angle and the addition formulas of the function *sin(l)* are well known, as follows:

$$\sin(2l) = 2\sin(l)\cos(l) \quad (or \quad sleaf_1(2l) = 2 sleaf_1(l) cleaf_1(l)) \tag{1.5}$$

$$\begin{aligned} &\sin(l_1 + l_2) = \sin(l_1)\cos(l_2) + \sin(l_2)\cos(l_1) \\ &(or \quad sleaf_1(l_1 + l_2) = sleaf_1(l_1) cleaf_1(l_2) + sleaf_1(l_2) cleaf_1(l_1)) \end{aligned} \tag{1.6}$$

For *n = 2*, Eq.(1.2) is rewritten as follows:

$$arcsleaf_2(r) = \int_0^r \frac{dt}{\sqrt{1-t^4}} = arcsl(r) = l \quad (1.7)$$

Function *arcsl(r)* represents the inverse lemniscatic elliptic function [12] [13] [14]. The double angle and the addition formulas are as follows:

$$sl(2l) = \frac{2sl(l)\sqrt{1-(sl(l))^4}}{1+(sl(l))^4} \quad \left( or \quad sleaf_2(2l) = \frac{2sleaf_2(l)\sqrt{1-(sleaf_2(l))^4}}{1+(sleaf_2(l))^4} \right)$$

(1.8)

$$sl(l_1 + l_2) = \frac{sl(l_1)\sqrt{1-(sl(l_2))^4} + sl(l_2)\sqrt{1-(sl(l_1))^4}}{1+(sl(l_1))^2(sl(l_2))^2}$$

$$\left( or \quad sleaf_2(l_1 + l_2) = \frac{sleaf_2(l_1)\sqrt{1-(sleaf_2(l_2))^4} + sleaf_2(l_2)\sqrt{1-(sleaf_2(l_1))^4}}{1+(sleaf_2(l_1))^2(sleaf_2(l_2))^2} \right)$$

(1.9)

In this study, the double angles with respect to the leaf functions *sleaf₃(l)* and *cleaf₃(l)* were derived by using the analogous of the trigonometric functions and lemniscatic elliptic function. Next, the addition formulas of these leaf functions were derived.

## 2 Double angle

### 2.1 Leaf function: *sleaf₃(l)*

For *n = 3*, the double angle is derived. Variables $l_1$ and $l_2$ are defined as follows:

$$l_1 = \int_0^{r_1} \frac{dt}{\sqrt{1-t^6}} \quad (2.1)$$

$$l_2 = \int_0^{r_2} \frac{dt}{\sqrt{1-t^6}} \quad (2.2)$$

Variable $l_1$ is related to variable $r_1$ as follows:

$$r_1 = sleaf_3(l_1) \quad (2.3)$$

Variable $l_2$ is related to variable $r_2$ as follows:

$$r_2 = sleaf_3(l_2) \tag{2.4}$$

We try to find the relationship between $r_1$ and $r_2$ when the following relationship exists.

$$l_1 = 2l_2 \tag{2.5}$$

The following equation is obtained by using Eqs. (2.1), (2.2), and (2.5).

$$\int_0^{r_1} \frac{dt}{\sqrt{1-t^6}} = 2\int_0^{r_2} \frac{dt}{\sqrt{1-t^6}} \tag{2.6}$$

The relationship between $r_1$ and $r_2$ is obtained as the double angle formula of the leaf function:

$$r_1 = \frac{2r_2\sqrt{1-r_2^6}}{\sqrt{1+8r_2^6}} \tag{2.7}$$

By differentiating it with respect to variable $r_2$, the following equation is obtained.

$$\frac{dr_1}{dr_2} = \frac{2 - 40r_2^6 - 16r_2^{12}}{(1+8r_2^6)^{\frac{3}{2}}\sqrt{1-r_2^6}} \tag{2.8}$$

The following equation is transformed by using the variable $r_2$.

$$\frac{1}{\sqrt{1-r_1^6}} = \pm \frac{(1+8r_2^6)^{\frac{3}{2}}}{1-20r_2^6-8r_2^{12}} \tag{2.9}$$

The following equation is obtained by using Eqs. (2.8) and (2.9).

$$\frac{1}{\sqrt{1-r_1^6}} \frac{dr_1}{dr_2} = \pm \frac{(1+8r_2^6)^{\frac{3}{2}}}{1-20r_2^6-8r_2^{12}} \frac{2-40r_2^6-16r_2^{12}}{(1+8r_2^6)^{\frac{3}{2}}\sqrt{1-r_2^6}} = \pm \frac{2}{\sqrt{1-r_2^6}} \tag{2.10}$$

The following equation is derived from the above equation.

$$\frac{dr_1}{\sqrt{1-r_1^6}} = 2\frac{dr_2}{\sqrt{1-r_2^6}} \tag{2.11}$$

Eq. (2.6) satisfies the above equation. The double angle of the leaf function is obtained as follows:

(1) In case of an inequality $\left(2m-\frac{1}{2}\right)\pi_3 \leq l \leq \left(2m+\frac{1}{2}\right)\pi_3$

(See [1] [2] for the constant $\pi_3$.)

$$sleaf_3(2l) = \frac{2 sleaf_3(l)\sqrt{1-(sleaf_3(l))^6}}{\sqrt{1+8(sleaf_3(l))^6}} \qquad (2.12)$$

(2) In case of an inequality $\left(2m+\frac{1}{2}\right)\pi_3 \leq l \leq \left(2m+\frac{3}{2}\right)\pi_3$

$$sleaf_3(2l) = -\frac{2 sleaf_3(l)\sqrt{1-(sleaf_3(l))^6}}{\sqrt{1+8(sleaf_3(l))^6}} \qquad (2.13)$$

The differentiation of the leaf function $sleaf_3(l)$ is defined by the domain of the variable $l$. The sign is defined as follows [1] [2].

(1) In case of the inequality $\left(2m-\frac{1}{2}\right)\pi_3 \leq l \leq \left(2m+\frac{1}{2}\right)\pi_3$

$$\frac{d}{dl} sleaf_3(l) = \sqrt{1-(sleaf_3(l))^6} \qquad (2.14)$$

(2) In case of the inequality $\left(2m+\frac{1}{2}\right)\pi_3 \leq l \leq \left(2m+\frac{3}{2}\right)\pi_3$

$$\frac{d}{dl} sleaf_3(l) = -\sqrt{1-(sleaf_3(l))^6} \qquad (2.15)$$

The sign of the term $\sqrt{1-(sleaf_3(l))^6}$ depends on the domain of variable $l$. Because of change in sign, it is necessary to classify Eqs. (2.12) and (2.13).

### 2.2 Leaf function: $cleaf_3(l)$

The leaf function $sleaf_3(l)$ is related to the leaf function $cleaf_3(l)$ as follows[1] [2]:

$$(sleaf_3(l))^2 + (cleaf_3(l))^2 + 2\cdot(sleaf_3(l))^2\cdot(cleaf_3(l))^2 = 1 \qquad (2.16)$$

Variable *l* in the above equation is replaced by variable *2l*.

$$(sleaf_3(2l))^2 + (cleaf_3(2l))^2 + 2 \cdot (sleaf_3(2l))^2 \cdot (cleaf_3(2l))^2 = 1 \tag{2.17}$$

By substituting Eq. (2.12) or Eq. (2.13) into Eq. (2.17), we get the following equation.

$$cleaf_3(2l) = \frac{2(cleaf_3(l))^2 + 2(cleaf_3(l))^4 - 1}{\sqrt{1 + 8(cleaf_3(l))^2 + 8(cleaf_3(l))^6 - 8(cleaf_3(l))^8}} \tag{2.18}$$

When the function *cleaf3(2l)* is derived from Eq. (2.17), both plus and minus signs occur in the function *cleaf3(2l)*. We have to choose the plus sign so that the initial condition, namely, *cleaf3(0)* =1 [1] [2] is satisfied.

## 3 Additional formulas

### 3.1 Leaf function: *sleaf3(l1+l2)*

With respect to arbitrary variables $l_1$ and $l_2$, the following equation is satisfied.

$$\{sleaf_3(l_1+l_2)\}^2 = \frac{\left\{sleaf_3(l_1)\frac{\partial sleaf_3(l_2)}{\partial l_2} + sleaf_3(l_2)\frac{\partial sleaf_3(l_1)}{\partial l_1}\right\}^2}{1 + 4(sleaf_3(l_2))^4(sleaf_3(l_1))^2 + 4(sleaf_3(l_2))^2(sleaf_3(l_1))^4} \\ + \frac{\{sleaf_3(l_2)(sleaf_3(l_1))^3 - (sleaf_3(l_2))^3 sleaf_3(l_1)\}^2}{1 + 4(sleaf_3(l_2))^4(sleaf_3(l_1))^2 + 4(sleaf_3(l_2))^2(sleaf_3(l_1))^4} \tag{3.1}$$

The above equation is set as follows:

$$\{g(l_1,l_2)\}^2 = \frac{\{p_1(l_1,l_2)\}^2}{p_3(l_1,l_2)} + \frac{\{p_2(l_1,l_2)\}^2}{p_3(l_1,l_2)} \tag{3.2}$$

$$g(l_1,l_2) = sleaf_3(l_1+l_2) \tag{3.3}$$

$$p_1(l_1,l_2) = sleaf_3(l_1)\frac{\partial sleaf_3(l_2)}{\partial l_2} + sleaf_3(l_2)\frac{\partial sleaf_3(l_1)}{\partial l_1} \tag{3.4}$$

$$p_2(l_1,l_2) = (sleaf_3(l_1))^3 sleaf_3(l_2) - sleaf_3(l_1)(sleaf_3(l_2))^3 \tag{3.5}$$

$$p_3(l_1,l_2) = 1 + 4(sleaf_3(l_1))^4(sleaf_3(l_2))^2 + 4(sleaf_3(l_1))^2(sleaf_3(l_2))^4 \tag{3.6}$$

For Eqs. (3.4)–(3.6), the following equations are obtained by differentiating with respect to variable $l_1$.

$$\frac{\partial p_1(l_1,l_2)}{\partial l_1} = \frac{\partial sleaf_3(l_1)}{\partial l_1}\frac{\partial sleaf_3(l_2)}{\partial l_2} - 3sleaf_3(l_2)(sleaf_3(l_1))^5 \tag{3.7}$$

$$\begin{aligned}\frac{\partial p_2(l_1,l_2)}{\partial l_1} &= 3(sleaf_3(l_1))^2 sleaf_3(l_2)\frac{\partial sleaf_3(l_1)}{\partial l_1} - (sleaf_3(l_2))^3\frac{\partial sleaf_3(l_1)}{\partial l_1} \\ &= sleaf_3(l_2)\{3(sleaf_3(l_1))^2 - (sleaf_3(l_2))^2\}\frac{\partial sleaf_3(l_1)}{\partial l_1}\end{aligned} \tag{3.8}$$

$$\begin{aligned}\frac{\partial p_3(l_1,l_2)}{\partial l_1} &= 16(sleaf_3(l_1))^3(sleaf_3(l_2))^2\frac{\partial sleaf_3(l_1)}{\partial l_1} + 8sleaf_3(l_1)(sleaf_3(l_2))^4\frac{\partial sleaf_3(l_1)}{\partial l_1} \\ &= 8sleaf_3(l_1)(sleaf_3(l_2))^2\{2(sleaf_3(l_1))^2 + (sleaf_3(l_2))^2\}\frac{\partial sleaf_3(l_1)}{\partial l_1}\end{aligned} \tag{3.9}$$

The following equation is substituted into Eq. (3.7).

$$\frac{\partial^2 sleaf_3(l)}{\partial l^2} = -3(sleaf_3(l))^5 \tag{3.10}$$

For Eqs. (3.4)–(3.6), the following equations are obtained by differentiating with respect to variable $l_2$.

$$\frac{\partial p_1(l_1,l_2)}{\partial l_2} = \frac{\partial sleaf_3(l_1)}{\partial l_1}\frac{\partial sleaf_3(l_2)}{\partial l_2} - 3sleaf_3(l_1)(sleaf_3(l_2))^5 \tag{3.11}$$

$$\begin{aligned}\frac{\partial p_2(l_1,l_2)}{\partial l_2} &= (sleaf_3(l_1))^3\frac{\partial sleaf_3(l_2)}{\partial l_2} - 3(sleaf_3(l_2))^2 sleaf_3(l_1)\frac{\partial sleaf_3(l_2)}{\partial l_2} \\ &= sleaf_3(l_1)\{(sleaf_3(l_1))^2 - 3(sleaf_3(l_2))^2\}\frac{\partial sleaf_3(l_2)}{\partial l_2}\end{aligned} \tag{3.12}$$

$$\begin{aligned}\frac{\partial p_3(l_1,l_2)}{\partial l_2} &= 16(sleaf_3(l_2))^3(sleaf_3(l_1))^2\frac{\partial sleaf_3(l_2)}{\partial l_2} + 8sleaf_3(l_2)(sleaf_3(l_1))^4\frac{\partial sleaf_3(l_2)}{\partial l_2} \\ &= 8sleaf_3(l_2)(sleaf_3(l_1))^2\{(sleaf_3(l_1))^2 + 2(sleaf_3(l_2))^2\}\frac{\partial sleaf_3(l_2)}{\partial l_2}\end{aligned} \tag{3.13}$$

For Eq. (3.2), the following equation is obtained by differentiating with respect to variable $l_1$.

$$\frac{\partial g(l_1,l_2)}{\partial l_1} = \frac{\left\{2p_1(l_1,l_2)\frac{\partial p_1(l_1,l_2)}{\partial l_1} + 2p_2(l_1,l_2)\frac{\partial p_2(l_1,l_2)}{\partial l_1}\right\}p_3(l_1,l_2)}{2g(l_1,l_2)(p_3(l_1,l_2))^2} \\ - \frac{\left\{(p_1(l_1,l_2))^2 + (p_2(l_1,l_2))^2\right\}\frac{\partial p_3(l_1,l_2)}{\partial l_1}}{2g(l_1,l_2)(p_3(l_1,l_2))^2}$$

(3.14)

Using Eqs. (3.7)–(3.9), the numerator of Eq. (3.14) is expanded as follows (see Appendix B):

$$\left\{2p_1(l_1,l_2)\frac{\partial p_1(l_1,l_2)}{\partial l_1} + 2p_2(l_1,l_2)\frac{\partial p_2(l_1,l_2)}{\partial l_1}\right\}p_3(l_1,l_2) \\ - \left\{(p_1(l_1,l_2))^2 + (p_2(l_1,l_2))^2\right\}\frac{\partial p_3(l_1,l_2)}{\partial l_1} \\ = \left\{2sleaf_3(l_1) - 8(sleaf_3(l_1))^5(sleaf_3(l_2))^2 - 24(sleaf_3(l_1))^3(sleaf_3(l_2))^4 \right. \\ \left. - 8sleaf_3(l_1)(sleaf_3(l_2))^6 - 16(sleaf_3(l_1))^5(sleaf_3(l_2))^8\right\}\frac{\partial sleaf_3(l_1)}{\partial l_1} \\ + \left\{2sleaf_3(l_2) - 8(sleaf_3(l_1))^2(sleaf_3(l_2))^5 - 24(sleaf_3(l_1))^4(sleaf_3(l_2))^3 \right. \\ \left. - 8(sleaf_3(l_1))^6 sleaf_3(l_2) - 16(sleaf_3(l_1))^8(sleaf_3(l_2))^5\right\}\frac{\partial sleaf_3(l_2)}{\partial l_2}$$

(3.15)

On the other hand, for Eq. (3.2), the following equation is obtained by differentiating with respect to variable $l_2$.

$$\frac{\partial g(l_1,l_2)}{\partial l_2} = \frac{\left\{2p_1(l_1,l_2)\frac{\partial p_1(l_1,l_2)}{\partial l_2} + 2p_2(l_1,l_2)\frac{\partial p_2(l_1,l_2)}{\partial l_2}\right\}p_3(l_1,l_2)}{2g(l_1,l_2)(p_3(l_1,l_2))^2} \\ - \frac{\left\{(p_1(l_1,l_2))^2 + (p_2(l_1,l_2))^2\right\}\frac{\partial p_3(l_1,l_2)}{\partial l_2}}{2g(l_1,l_2)(p_3(l_1,l_2))^2}$$

(3.16)

Using Eqs. (3.11)–(3.13), the numerator of Eq. (3.16) is expanded as follows:

$$\left\{2p_1(l_1,l_2)\frac{\partial p_1(l_1,l_2)}{\partial l_2}+2p_2(l_1,l_2)\frac{\partial p_2(l_1,l_2)}{\partial l_2}\right\}p_3(l_1,l_2)-\left\{(p_1(l_1,l_2))^2+(p_2(l_1,l_2))^2\right\}\frac{\partial p_3(l_1,l_2)}{\partial l_2}$$

$$=\left\{2sleaf_3(l_1)-8(sleaf_3(l_1))^5(sleaf_3(l_2))^2-24(sleaf_3(l_1))^3(sleaf_3(l_2))^4\right.$$

$$\left.-8sleaf_3(l_1)(sleaf_3(l_2))^6-16(sleaf_3(l_1))^5(sleaf_3(l_2))^8\right\}\frac{\partial sleaf_3(l_1)}{\partial l_1}$$

$$+\left\{2sleaf_3(l_2)-8(sleaf_3(l_1))^2(sleaf_3(l_2))^5-24(sleaf_3(l_1))^4(sleaf_3(l_2))^3\right.$$

$$\left.-8(sleaf_3(l_1))^6 sleaf_3(l_2)-16(sleaf_3(l_1))^8(sleaf_3(l_2))^5\right\}\frac{\partial sleaf_3(l_2)}{\partial l_2}$$

(3.17)

Using the above equation, the following equation is derived.

$$\frac{\partial g(l_1,l_2)}{\partial l_1}=\frac{\partial g(l_1,l_2)}{\partial l_2} \tag{3.18}$$

The following equation that satisfies Eq. (3.18) (see Appendix A) is derived.

$$g(l_1,l_2)=g(l_1+l_2,0) \tag{3.19}$$

Using the initial condition $sleaf_3(0)=0$ and $\partial sleaf_3(0)/\partial l=1$, the function $g(l_1+l_2, 0)$ is written as follows:

$$\{g(l_1+l_2,0)\}^2=\frac{\left\{sleaf_3(l_1+l_2)\frac{\partial sleaf_3(0)}{\partial l_2}+sleaf_3(0)\frac{\partial sleaf_3(l_1+l_2)}{\partial l_1}\right\}^2}{1+4(sleaf_3(l_1+l_2))^4(sleaf_3((0)))^2+4(sleaf_3(l_1+l_2))^2(sleaf_3(0))^4}$$

$$+\frac{\left\{(sleaf_3(l_1+l_2))^3 sleaf_3(0)-sleaf_3(l_1+l_2)(sleaf_3(0))^3\right\}^2}{1+4(sleaf_3(l_1+l_2))^4(sleaf_3((0)))^2+4(sleaf_3(l_1+l_2))^2(sleaf_3(0))^4}$$

$$=\{sleaf_3(l_1+l_2)\}^2$$

(3.20)

Using Eqs. (3.19) and (3.20), the following equation is obtained.

$$\{sleaf_3(l_1+l_2)\}^2 = g(l_1+l_2, 0) = g(l_1, l_2)$$

$$= \frac{\left\{sleaf_3(l_1) \cdot \frac{\partial sleaf_3(l_2)}{\partial l_2} + sleaf_3(l_2) \cdot \frac{\partial sleaf_3(l_1)}{\partial l_1}\right\}^2}{1 + 4(sleaf_3(l_1))^4(sleaf_3(l_2))^2 + 4(sleaf_3(l_1))^2(sleaf_3(l_2))^4}$$

$$+ \frac{\{(sleaf_3(l_1))^3 sleaf_3(l_2) - sleaf_3(l_1)(sleaf_3(l_2))^3\}^2}{1 + 4(sleaf_3(l_1))^4(sleaf_3(l_2))^2 + 4(sleaf_3(l_1))^2(sleaf_3(l_2))^4}$$

(3.21)

The differentiation $\partial sleaf_3(l_1)/dl_1$ and $\partial sleaf_3(l_2)/dl_2$ depends on the domains of variables $l_1$ and $l_2$. It is necessary to classify domains. The additional formulas can be summarized as follows.

(i) In case both $(4m-1)\pi_3/2 \leqq l_1 \leqq (4m+1)\pi_3/2$ and $(4k-1)\pi_3/2 \leqq l_2 \leqq (4k+1)\pi_3/2$ or both $(4m+1)\pi_3/2 \leqq l_1 \leqq (4m+3)\pi_3/2$ and $(4k+1)\pi_3/2 \leqq l_2 \leqq (4k+3)\pi_3/2$ (Variables $m$ and $k$ represent integers; for the constant $\pi_3$, see [1][2].)

$$\{sleaf_3(l_1+l_2)\}^2 = \frac{\left\{sleaf_3(l_1)\sqrt{1-(sleaf_3(l_2))^6} + sleaf_3(l_2)\sqrt{1-(sleaf_3(l_1))^6}\right\}^2}{1 + 4(sleaf_3(l_1))^4(sleaf_3(l_2))^2 + 4(sleaf_3(l_1))^2(sleaf_3(l_2))^4}$$

$$+ \frac{\{(sleaf_3(l_1))^3 sleaf_3(l_2) - sleaf_3(l_1)(sleaf_3(l_2))^3\}^2}{1 + 4(sleaf_3(l_1))^4(sleaf_3(l_2))^2 + 4(sleaf_3(l_1))^2(sleaf_3(l_2))^4}$$

(3.22)

(ii) In case both $(4m+1)\pi_3/2 \leqq l_1 \leqq (4m+3)\pi_3/2$ and $(4k-1)\pi_3/2 \leqq l_2 \leqq (4k+1)\pi_3/2$ or both $(4m-1)\pi_3/2 \leqq l_1 \leqq (4m+1)\pi_3/2$ and $(4k+1)\pi_3/2 \leqq l_2 \leqq (4k+3)\pi_3/2$ (Variables $m$ and $k$ represent integers)

$$\{sleaf_3(l_1+l_2)\}^2 = \frac{\left\{sleaf_3(l_1)\sqrt{1-(sleaf_3(l_2))^6} - sleaf_3(l_2)\sqrt{1-(sleaf_3(l_1))^6}\right\}^2}{1 + 4(sleaf_3(l_1))^4(sleaf_3(l_2))^2 + 4(sleaf_3(l_1))^2(sleaf_3(l_2))^4}$$

$$+ \frac{\{(sleaf_3(l_1))^3 sleaf_3(l_2) - sleaf_3(l_1)(sleaf_3(l_2))^3\}^2}{1 + 4(sleaf_3(l_1))^4(sleaf_3(l_2))^2 + 4(sleaf_3(l_1))^2(sleaf_3(l_2))^4}$$

(3.23)

### 3.2 Leaf function: *cleaf₃(l₁+l₂)*

With respect to arbitrary variables $l_1$ and $l_2$, the following equation is satisfied.

$$(cleaf_3(l_1+l_2))^2 = \frac{\left\{cleaf_3(l_1)\frac{\partial sleaf_3(l_2)}{\partial l_2}+sleaf_3(l_2)\frac{\partial cleaf_3(l_1)}{\partial l_1}\right\}^2}{1+4(sleaf_3(l_2))^4(cleaf_3(l_1))^2+4(sleaf_3(l_2))^2(cleaf_3(l_1))^4}$$
$$+\frac{\left\{(sleaf_3(l_1))^3 cleaf_3(l_2)-sleaf_3(l_1)(cleaf_3(l_2))^3\right\}^2}{1+4(sleaf_3(l_2))^4(cleaf_3(l_1))^2+4(sleaf_3(l_2))^2(cleaf_3(l_1))^4} \quad (3.24)$$

The additional formulas of the leaf function *cleaf₃(l₁ + l₂)* can be similarly proved as described in section 3.1. The differentials ∂*cleaf₃(l₁)/dl₁* and ∂*cleaf₃(l₂)/dl₂* depend on the domains of variables *l₁* and *l₂*. It is necessary to classify domains. The additional formulas can be summarized as follows.

(i) In case both *2kπ₃ ≦ l₁ ≦ (2k+1)π₃* and *(4m-1)π₃/2 ≦ l₂ ≦ (4m+1)π₃/2* or *(2k+1)π₃ ≦ l₁ ≦ (2k+2)π₃* and *(4m+1)π₃/2 ≦ l₂ ≦ (4m+3)π₃/2* (Variables *m* and *k* represent integers)

$$(cleaf_3(l_1+l_2))^2 = \frac{\left\{cleaf_3(l_1)\sqrt{1-(sleaf_3(l_2))^6}-sleaf_3(l_2)\sqrt{1-(cleaf_3(l_1))^6}\right\}^2}{1+4(sleaf_3(l_2))^4(cleaf_3(l_1))^2+4(sleaf_3(l_2))^2(cleaf_3(l_1))^4}$$
$$+\frac{\left\{(sleaf_3(l_1))^3 cleaf_3(l_2)-sleaf_3(l_2)(cleaf_3(l_1))^3\right\}^2}{1+4(sleaf_3(l_2))^4(cleaf_3(l_1))^2+4(sleaf_3(l_2))^2(cleaf_3(l_1))^4} \quad (3.25)$$

(ii) In case both *(2k+1)π₃ ≦ l₁ ≦ (2k+2)π₃* and *(4m-1)π₃/2 ≦ l₂ ≦ (4m+1)π₃/2* or *2kπ₃ ≦ l₁ ≦ (2k+1)π₃* and *(4m+1)π₃/2 ≦ l₂ ≦ (4m+3)π₃/2* (Variables *m* and *k* represent integers)

$$(cleaf_3(l_1+l_2))^2 = \frac{\left\{cleaf_3(l_1)\sqrt{1-(sleaf_3(l_2))^6}+sleaf_3(l_2)\sqrt{1-(cleaf_3(l_1))^6}\right\}^2}{1+4(sleaf_3(l_2))^4(cleaf_3(l_1))^2+4(sleaf_3(l_2))^2(cleaf_3(l_1))^4}$$
$$+\frac{\left\{(sleaf_3(l_1))^3 cleaf_3(l_2)-sleaf_3(l_2)(cleaf_3(l_1))^3\right\}^2}{1+4(sleaf_3(l_2))^4(cleaf_3(l_1))^2+4(sleaf_3(l_2))^2(cleaf_3(l_1))^4} \quad (3.26)$$

## 4 Numerical example

The numerical data of the leaf functions *sleaf₃(l)* and *cleaf₃(l)* are presented in Table 1. Using Eqs. (3.22), (3.23), (3.25), and (3.26), the values of *sleaf₃(l₁+ l₂)* and *cleaf₃(l₁+l₂)* can be calculated by using the values of *sleaf₃(l₁)*, *sleaf₃(l₂)*, *cleaf₃(l₁)*, and *cleaf₃(l₂)*. As an example, consider *l₁ = 0.2* and *l₂ = 0.3*.

Then, *sleaf₃(0.2+0.3)* can be obtained. We have to choose Eq. (3.22) for the inequalities $-\pi_3/2 \leqq 0.2 \leqq \pi_3/2$ and $-\pi_3/2 \leqq 0.3 \leqq \pi_3/2$. By substituting $l_1 = 0.2$ and $l_2 = 0.3$ into Eq. (3.22), we get the following equation:

$$\{sleaf_3(0.2+0.3)\}^2 = \frac{\{sleaf_3(0.2)\sqrt{1-(sleaf_3(0.3))^6} + sleaf_3(0.3)\sqrt{1-(sleaf_3(0.2))^6}\}^2}{1+4(sleaf_3(0.2))^4(sleaf_3(0.3))^2 + 4(sleaf_3(0.2))^2(sleaf_3(0.3))^4}$$

$$+ \frac{\{(sleaf_3(0.2))^3 sleaf_3(0.3) - sleaf_3(0.2)(sleaf_3(0.3))^3\}^2}{1+4(sleaf_3(0.2))^4(sleaf_3(0.3))^2 + 4(sleaf_3(0.2))^2(sleaf_3(0.3))^4}$$

$$= \frac{\{(0.199999\cdots)\sqrt{1-(0.299984\cdots)^6} + (0.299984\cdots)\sqrt{1-(0.199999\cdots)^6}\}^2}{1+4(0.199999\cdots)^4(0.299984\cdots)^2 + 4(0.199999\cdots)^2(0.299984\cdots)^4}$$

$$+ \frac{\{(0.199999\cdots)^3(0.299984\cdots) - (0.199999\cdots)(0.299984\cdots)^3\}^2}{1+4(0.199999\cdots)^4(0.299984\cdots)^2 + 4(0.199999\cdots)^2(0.299984\cdots)^4} = 0.2494431\cdots$$

(4.1)

Therefore, the value of *sleaf₃(0.5)* is obtained as follows:

$$sleaf_3(0.2+0.3) = sleaf_3(0.5) = \sqrt{0.2494431\cdots} = 0.49944\cdots \qquad (4.2)$$

Table 1 Numerical data of leaf functions *sleaf₃(l)* and *cleaf₃(l)*

| l | sleaf₃(l) | cleaf₃(l) | l | sleaf₃(l) | cleaf₃(l) |
|---|---|---|---|---|---|
| 0.0 | 0.000000 | 1.000000 | 2.1 | 0.328621 | -0.856490 |
| 0.1 | 0.100000 | 0.985184 | 2.2 | 0.228649 | -0.926290 |
| 0.2 | 0.199999 | 0.942810 | 2.3 | 0.128651 | -0.975670 |
| 0.3 | 0.299984 | 0.878184 | 2.4 | 0.028651 | -0.998770 |
| 0.4 | 0.399883 | 0.797825 | 2.5 | -0.071350 | -0.992410 |
| 0.5 | 0.499443 | 0.707632 | 2.6 | -0.171350 | -0.957500 |
| 0.6 | 0.598009 | 0.611979 | 2.7 | -0.271340 | -0.898590 |
| 0.7 | 0.694183 | 0.513647 | 2.8 | -0.371280 | -0.822090 |
| 0.8 | 0.785387 | 0.414176 | 2.9 | -0.470980 | -0.734190 |
| 0.9 | 0.867486 | 0.314304 | 3.0 | -0.569930 | -0.639750 |
| 1.0 | 0.934768 | 0.214324 | 3.1 | -0.667000 | -0.541980 |
| 1.1 | 0.980708 | 0.114325 | 3.2 | -0.759970 | -0.442740 |
| 1.2 | 0.999692 | 0.014325 | 3.3 | -0.845200 | -0.342940 |
| 1.3 | 0.989090 | -0.085670 | 3.4 | -0.917390 | -0.242970 |
| 1.4 | 0.950393 | -0.185670 | 3.5 | -0.970090 | -0.142980 |
| 1.5 | 0.888560 | -0.285660 | 3.6 | -0.997240 | -0.042980 |
| 1.6 | 0.810064 | -0.385580 | 3.7 | -0.995140 | 0.057024 |
| 1.7 | 0.720972 | -0.485220 | 3.8 | -0.964110 | 0.157024 |
| 1.8 | 0.625896 | -0.583990 | 3.9 | -0.908270 | 0.257019 |
| 1.9 | 0.527828 | -0.680640 | 4.0 | -0.833880 | 0.356971 |
| 2.0 | 0.428461 | -0.772770 | 4.1 | -0.747280 | 0.456727 |

# 5 Conclusions

By using the analogy of the inverse trigonometric and lemniscate functions, higher order of functions were artificially created as leaf functions. The following conclusions can be drawn from this study:

・The artificially created leaf function $r = sleaf_3(l)$ corresponds to $r = sin(l)$ and $r = sl(l)$, whereas the artificially created leaf function $r = cleaf_3(l)$ corresponds to $r = cos(l)$ and $r = cl(l)$. The waves obtained from the trigonometric functions $sin(l)$ and $cos(l)$ have a periodicity of $6.28\cdots$ and an amplitude of 1, while the waves obtained from the lemniscate functions have a periodicity of $5.24\cdots$ and an amplitude of 1. The artificially created functions $sleaf_3(l)$ and $cleaf_3(l)$ also produce regular waves with $4.85\cdots$ periodicity and amplitude of 1. These functions produce waves different from those of both the trigonometric and lemniscate functions. The prefix "s" and "c" of these functions are classified by the following initial conditions: $r(0) = 0$, $dr(0)/dl = 1$ and $r(0) = 1$, $dr(0)/dl = 0$.

・Additional formulas of leaf functions $sleaf_3(l_1 + l_2)$ and $cleaf_3(l_1 + l_2)$ were also derived.

・The values of the leaf functions $sleaf_3(l_1 + l_2)$ can be numerically obtained using the values of $sleaf_3(l_1)$, $sleaf_3(l_2)$, while those of $cleaf_3(l_1 + l_2)$ can be obtained using values of $sleaf_3(l_1)$, $sleaf_3(l_2)$, $cleaf_3(l_1)$, and $cleaf_3(l_2)$ through the additional formulas.

This study will be the first step in clarifying solution behaviors of nonlinear equations consisting of both the second-order differentials and the higher order polynomials.

## Appendix A

Let $g(l_1, l_2)$ be a differentiable function. The necessary and sufficient condition for satisfying $g(l_1, l_2) = g(l_1+l_2, 0)$ is that $\frac{\partial g(l_1, l_2)}{\partial l_1} = \frac{\partial g(l_1, l_2)}{\partial l_2}$ holds. Function $h(x,y)$ is defined as follows:

$$h(x, y) = g(x + y, x - y) \tag{A.1}$$

$$l_1 = x + y \tag{A.2}$$

$$l_2 = x - y \tag{A.3}$$

By differentiating the above equation with respect to $y$, we obtain the following equation.

$$\begin{aligned}\frac{\partial h(x,y)}{\partial y} &= \frac{\partial g(x+y, x-y)}{\partial l_1}\frac{\partial l_1}{\partial y} + \frac{\partial g(x+y, x-y)}{\partial l_2}\frac{\partial l_2}{\partial y} \\ &= \frac{\partial g(x+y, x-y)}{\partial l_1} - \frac{\partial g(x+y, x-y)}{\partial l_2}\end{aligned} \tag{A.4}$$

Therefore, if the equation $\partial g/\partial l_1 = \partial g/\partial l_2$ holds, the following equation holds.

$$\frac{\partial h(x,y)}{\partial y} = 0 \tag{A.5}$$

From the above equation, we find that $h(x,y)$ is a function of $x$, not a function of $y$. Therefore, the following equation holds for any constants $a$ and $b$.

$$h(x,a) = h(x,b) \tag{A.6}$$

Here, we set the following equation:

$$x = b = \frac{l_1 + l_2}{2} \tag{A.7}$$

$$a = \frac{l_1 - l_2}{2} \tag{A.8}$$

From the above two equations, the following equation is obtained.

$$h(x,a) = h\left(\frac{l_1+l_2}{2}, \frac{l_1-l_2}{2}\right) = g\left(\frac{l_1+l_2}{2} + \frac{l_1-l_2}{2}, \frac{l_1+l_2}{2} - \frac{l_1-l_2}{2}\right) = g(l_1, l_2) \tag{A.9}$$

$$h(x,b) = h\left(\frac{l_1+l_2}{2}, \frac{l_1+l_2}{2}\right) = g\left(\frac{l_1+l_2}{2} + \frac{l_1+l_2}{2}, \frac{l_1+l_2}{2} - \frac{l_1+l_2}{2}\right) = g(l_1+l_2, 0) \tag{A.10}$$

The following equation is obtained by using Eqs. (A.6), (A.9), and (A.10).

$$g(l_1, l_2) = g(l_1 + l_2, 0) \tag{A.11}$$

Conversely, if the above equation holds, the following relational expression can be obtained by using Eqs. (A.1) and (A.11).

$$h(x, y) = g(x + y, x - y) = g(2x, 0) \tag{A.12}$$

We differentiate the above equation with $y$. Because $g(2x, 0)$ is a function of $x$ and not a function of $y$, the following equation is obtained.

$$\frac{\partial}{\partial y} h(x, y) = \frac{\partial}{\partial y} g(2x, 0) = 0 \tag{A.13}$$

$$\begin{aligned}\frac{\partial h(x, y)}{\partial y} &= \frac{\partial g(x+y, x-y)}{\partial l_1} \frac{\partial l_1}{\partial y} + \frac{\partial g(x+y, x-y)}{\partial l_2} \frac{\partial l_2}{\partial y} \\ &= \frac{\partial g(x+y, x-y)}{\partial l_1} - \frac{\partial g(x+y, x-y)}{\partial l_2} = 0\end{aligned} \tag{A.14}$$

Therefore, the following equation is obtained.

$$\frac{\partial g(x+y, x-y)}{\partial l_1} = \frac{\partial g(x+y, x-y)}{\partial l_2} \tag{A.15}$$

Here, we use the following equation:

$$x = \frac{l_1 + l_2}{2} \tag{A.16}$$

$$y = \frac{l_1 - l_2}{2} \tag{A.17}$$

Using Eq. (A.15), we obtain the following equation:

$$\frac{\partial g(l_1, l_2)}{\partial l_1} = \frac{\partial g(l_1, l_2)}{\partial l_2} \tag{A.18}$$

## Appendix B

The first term of the Eq. (3.15) is transformed as follows:

$$
\begin{aligned}
& p_1(l_1,l_2)\frac{\partial p_1(l_1,l_2)}{\partial l_1} + p_2(l_1,l_2)\frac{\partial p_2(l_1,l_2)}{\partial l_1} \\
&= sleaf_3(l_1)\frac{\partial sleaf_3(l_1)}{\partial l_1}\left(\frac{\partial sleaf_3(l_2)}{\partial l_2}\right)^2 + sleaf_3(l_2)\left(\frac{\partial sleaf_3(l_1)}{\partial l_1}\right)^2\frac{\partial sleaf_3(l_2)}{\partial l_2} \\
&\quad - 3(sleaf_3(l_1))^6 sleaf_3(l_2)\frac{\partial sleaf_3(l_2)}{\partial l_2} - 4(sleaf_3(l_1))^3(sleaf_3(l_2))^4\frac{\partial sleaf_3(l_1)}{\partial l_1} \\
&\quad + sleaf_3(l_1)(sleaf_3(l_2))^6\frac{\partial sleaf_3(l_1)}{\partial l_1} \\
&= sleaf_3(l_1)\frac{\partial sleaf_3(l_1)}{\partial l_1}\left(1-(sleaf_3(l_2))^6\right) + sleaf_3(l_2)\left(1-(sleaf_3(l_1))^6\right)\frac{\partial sleaf_3(l_2)}{\partial l_2} \\
&\quad - 3(sleaf_3(l_1))^6 sleaf_3(l_2)\frac{\partial sleaf_3(l_2)}{\partial l_2} - 4(sleaf_3(l_2))^4(sleaf_3(l_1))^3\frac{\partial sleaf_3(l_1)}{\partial l_1} \\
&\quad + sleaf_3(l_1)(sleaf_3(l_2))^6\frac{\partial sleaf_3(l_1)}{\partial l_1} \\
&= \left\{sleaf_3(l_1) - 4(sleaf_3(l_1))^3(sleaf_3(l_2))^4\right\}\frac{\partial sleaf_3(l_1)}{\partial l_1} \\
&\quad + \left\{sleaf_3(l_2) - 4(sleaf_3(l_1))^6 sleaf_3(l_2)\right\}\frac{\partial sleaf_3(l_2)}{\partial l_2}
\end{aligned}
\tag{B.1}
$$

The following equation is applied [1] [2].

$$
\frac{d}{dl}sleaf_3(l) = \sqrt{1-(sleaf_3(l))^6}
\tag{B.2}
$$